\documentclass{article}
\usepackage[francais]{babel}
\usepackage[latin1]{inputenc}
\usepackage[T1]{fontenc}
\usepackage{amsmath}
\usepackage[psamsfonts]{amsfonts}
\usepackage{amsthm}
\usepackage{amsopn}
\usepackage{amscd}
\usepackage{dsfont}
\usepackage{makeidx}
\usepackage[dvipdf]{graphics}
\usepackage{psfrag}
\usepackage{pict2e}
\usepackage[left=5cm, right=5cm, top=3cm, bottom=4cm]{geometry}
\newcommand{\Cplx}{\mathbb C}

\theoremstyle{plain}
\newtheorem{Theo}{Théorème}

\newtheorem*{Lem}{Lemme}
\newtheorem*{Lemgeom}{Lemme géométrique}

\theoremstyle{definition}
\newtheorem{Cor}{Corollaire}
\newtheorem{Prop}[Theo]{Proposition}

\newtheorem*{Fait}{Fait 1}
\newtheorem*{Faitt}{Fait 2}
\newtheorem*{Faittt}{Fait 3}
\newtheorem*{Faitttt}{Fait 4}
\theoremstyle{remark}
\newtheorem*{demo}{démonstration}

\catcode`\«=\active \catcode`\»=\active
\def«{\og\ignorespaces}
\def»{{\fg}}
\begin{document}

\begin{center}
{\LARGE\bf Un théorème de Bloch presque complexe}\\
\vspace{0.5cm}
\large Benoît Saleur\footnote{Département de Mathématiques de la faculté des sciences d'Orsay, Université Paris-Sud 11, 91405 Orsay Cedex\\
benoit.saleur@math.u-psud.fr\\
Mots clefs: hyperbolicité complexe, théorie de Nevanlinna, Bloch, courbes pseudoholomorphes, courants positifs\\
Clefs AMS: 32A18, 32H30, 32Q45, 32Q65, 32U40}
\end{center}

\begin{center}
{{\large \bf Résumé.}
Etant donné un plan projectif presque complexe $(\mathbb{P}^{2}(\Cplx ),J)$, nous démontrons que toute suite non normale de $J$-disques évitant une configuration $C$ de quatre $J$-droites en position générale admet une sous-suite convergeant, au sens de Hausdorff, vers la réunion $\Delta$ de trois $J$-droites. En particulier, le complémentaire de quatre $J$-droites en position générale est hyperboliquement plongé dans $(\mathbb{P}^{2}(\Cplx ),J)$ modulo $\Delta$.}
\end{center}

\section{Introduction}
Soit $X$ une variété de dimension $2n$ munie d'une structure presque complexe, c'est-à dire d'un automorphisme lisse $J$ de $TX$ tel que $J^{2}=-Id$. Même lorsque $J$ n'est pas intégrable, la variété $X$ possède de nombreuses courbes $J$-holomorphes, i.e. des surfaces dont le plan tangent est en tout point une droite complexe pour $J$. Le théorème d'uniformisation des surfaces de Riemann assure qu'une telle courbe est paramétrée localement par un $J$-disque, autrement dit par une application $f:(\mathbb{D},i)\rightarrow (X,J)$ définie sur le disque unité de $\Cplx$ et $J$-holomorphe, i.e. vérifiant $df\circ i = J\circ df$.\\
Rappelons qu'un $J$-disque est de classe $\mathcal{C}^{\infty}$, et que la limite uniforme locale d'une suite de $J$-disques est encore un $J$-disque (voir l'article de J.-C. Sikorav dans \cite{Audin}). Une suite $(f_{n})$ de $J$-disques est dite normale si elle admet un sous-suite convergeant localement uniformément vers un $J$-disque. \\

L'existence, pour tout point $P\in X$ et tout vecteur $v\in T_{P}X$, d'un $J$-disque $f$ non constant passant par $P$ tangentiellement à $v$ (voir l'article de J.-C. Sikorav dans \cite{Audin}), motive la définition d'une pseudométrique de Kobayashi-Royden $K_{X}$(voir \cite{Kruglikov}):

$$K_{X}(P,v)=\inf \left\{ \frac{1}{|\lambda |} >0 \text { | il existe } f:\mathbb{D} \rightarrow X \hspace{0.3cm}J\text{-holomorphe avec } f(0)=P, df_{0}(\frac{\partial}{\partial z})=\lambda v \right\}.$$

La variété $X$ est dite hyperbolique au sens de Kobayashi lorsque $K_{X}$ est non dégénérée. Dans le cas contraire, il existe un point par lequel passent des $J$-disques arbitrairement grands dans une direction donnée. Si la variété $X$ est compacte, cela équivaut à l'existence d'une courbe entière, c'est à dire d'une application $J$-holomorphe $f:(\Cplx ,i)\rightarrow (X,J)$ non constante. Ce critère, dû à Brody (voir \cite{Brody}), découle du théorème de reparamétrisation suivant:

\begin{Theo}\label{Brody}[Théorème de reparamétrisation de Brody]---
Soit $f_{n}:\mathbb{D}\rightarrow X$ une suite non normale de $J$-disques. Il existe une suite de contractions affines $\rho_{n}$ de $\Cplx$ convergeant vers un point du disque unité, appelé point d'explosion, telle que $f_{n}\circ \rho_{n}$ converge uniformément sur tout compact de $\Cplx$ vers une courbe de Brody, c'est-à-dire une courbe entière non constante à dérivée bornée.
\end{Theo}

Le critère de Brody reste valable lorsque $X$ est le complémentaire dans l'espace projectif complexe $(\mathbb{P}^{n}(\Cplx ),i)$ de $2n+1$ hyperplans en position générale: en vertu du théorème de Green interdisant l'existence d'une courbe entière dans le domaine $X$ (voir \cite{Green}), la pseudométrique $K_{X}$ est donc non dégénérée, et même supérieure à la restriction à $X$ d'une certaine métrique définie sur $\mathbb{P}^{n}(\Cplx )$ (le domaine $X$ est dit hyperboliquement plongé dans l'espace projectif).\\
Le complémentaire de $p<2n+1$ hyperplans de $\mathbb{P}^{n}(\Cplx )$ contient des courbes entières, et n'est donc pas hyperbolique au sens de Kobayashi. Cependant, A. Bloch puis H. Cartan ont montré que la suite des images $f_{n}(\mathbb{D})$ d'une suite non normale de disques holomorphes $f_{n}:\mathbb{D}\rightarrow \mathbb{P}^{n}(\Cplx )$ évitant $n+2$ hyperplans $H_{k}$ ($0\leq k \leq n+1$) en position générale converge au sens de Hausdorff vers une union finie d'hyperplans dits diagonaux (en identifiant $\mathbb{P}^{n}(\Cplx )$ avec l'hyperplan de $\mathbb{P}^{n+1}(\Cplx )$ d'équation $X_{1} + ... + X_{n+1}=0$ et les $H_{k}$ avec les axes de coordonnées, les hyperplans diagonaux sont donnés par les équations $\sum_{l\in I}X_{l}=0$, pour les parties $I$ de $\{0,...,n+1 \}$ de cardinal $2\leq \# I \leq n$). En particulier la pseudométrique de Kobayashi du complémentaire des  $H_{k}$ dégénère exactement sur les hyperplans diagonaux, et toute courbe entière non constante est portée par un hyperplan diagonal. Cette dernière propriété est connue sous le nom de théorème de Borel.\\

Une version presque complexe du théorème de Green est connue depuis J. Duval en dimension réelle $4$ (voir \cite{Green presque complexe}), la question ne se posant pas en dimension supérieure, à cause du manque d'hypersurfaces $J$-holomorphes. Rappelons brièvement son résultat.\\
Soit $J$ une structure presque complexe de classe $\mathcal{C}^{\infty}$ sur $\mathbb{P}^{2}(\Cplx )$ positive par rapport à la métrique de Fubini-Study $\omega$, i.e. $\omega_{P} (v,Jv)>0$ pour tout point $P$ et tout vecteur $v\in T_{P}\mathbb{P}^{2}(\Cplx )$ non nul. Une $J$-droite du plan projectif presque complexe $(\mathbb{P}^{2}(\Cplx ),J)$ est une courbe $J$-holomorphe plongée dans $\mathbb{P}^{2}(\Cplx )$, difféomorphe à $\mathbb{P}^{1}(\Cplx )$ et de degré $1$ en homologie. D'après M. Gromov \cite{Gromov} (voir aussi J.-C. Sikorav \cite{Sikorav}), l'espace des $J$-droites est difféomorphe à $\mathbb{P}^{2}(\Cplx )$. En outre, par deux points distincts passe une unique $J$-droite, deux $J$-droites distinctes se coupent transversalement en un unique point, et les $J$-droites passant par un point $P$ forment un pinceau difféomorphe à $\mathbb{P}^{1}(\Cplx )$, donnant une projection centrale $\pi_{P}:\mathbb{P}^{2}(\Cplx )\setminus \{P\} \rightarrow \mathbb{P}^{1}(\Cplx )$.\\
Dans ce contexte le théorème de Green s'énonce ainsi:

\begin{Theo}\cite{Green presque complexe}---
Le complémentaire dans $(\mathbb{P}^{2}(\Cplx ),J)$ d'une réunion de cinq $J$-droites en position générale (i.e. sans point triple) est hyperbolique au sens de Kobayashi.
\end{Theo}

Il est naturel de poursuivre l'étude et de s'intéresser au complémentaire dans $(\mathbb{P}^{2}(\Cplx ),J)$ de quatre $J$-droites $L_{1},L_{2}, L_{3}$ et $L_{4}$ en position générale. Nous noterons $S_{k}$ ($1\leq k \leq 6$) les six points doubles de la configuration $C=L_{1}\cup L_{2} \cup L_{3} \cup L_{4}$. Appelons diagonales les trois $J$-droites  $\Delta_{1}$, $\Delta_{2}$ et $\Delta_{3}$ coupant chacune $C$ en deux points doubles. Leur réunion $\Delta$ est appelée diviseur diagonal.\\
Voici l'énoncé que nous allons démontrer:

\begin{Theo}\label{Bloch}
Soit $f_{n}:\mathbb{D}\rightarrow \mathbb{P}^{2}(\Cplx )\setminus C$ une suite non normale de $J$-disques évitant la configuration $C$. Alors $(f_{n}(\mathbb{D}))$ converge au sens de Hausdorff vers le diviseur diagonal.
\end{Theo}

Comme dans le cas complexe, la pseudométrique $K_{\mathbb{P}^{2}(\Cplx )\setminus C}$ est uniformément nulle tangentiellement aux diagonales (puisque celles-ci ne coupent $C$ qu'en deux points, elles contiennent des courbes entières), mais ne peut dégénérer hors de $\Delta$. Le complémentaire de $C$ est dit hyperboliquement plongé dans le plan projectif presque complexe modulo le diviseur diagonal. Ceci implique en particulier un théorème de Borel presque complexe:

\begin{Cor}
Soit $f:\Cplx \rightarrow \mathbb{P}^{2}(\Cplx )\setminus C$ une courbe entière non constante évitant la configuration $C$. Alors $f(\Cplx )$ est contenue dans l'une des diagonales.
\end{Cor}

En effet, si $f:\Cplx \rightarrow \mathbb{P}^{2}(\Cplx )\setminus C$ est une courbe entière non constante, pour tout point $z$ de $\Cplx$, $K_{\mathbb{P}^{2}(\Cplx )\setminus C}(f(z), df_{z}(v))\leq K_{\Cplx}(z,v)=0$, donc $f(z)\in \Delta$ dès lors que $df_{z}$ est non nulle. La courbe $f(\Cplx )$ est donc contenue dans l'une des diagonales.\\

Remarquons que le théorème de Bloch ne se réduit pas à ce dernier résultat: les courbes entières issues des explosions d'une suite $(f_{n})$ peuvent être toutes contenues dans une $J$-droite $L$ sans pour autant que la suite des $f_{n}(\mathbb{D})$ converge au sens de Hausdorff vers $L$. Ce défaut de localisation inhérent au théorème de reparamétrisation de Brody nous imposera des précautions particulières.\\ 

Nous allons démontrer le théorème \ref{Bloch} par l'absurde: supposant donnée une suite non normale de $J$-disques évitant $C$ et ne convergeant pas vers les diagonales, nous construirons un courant positif fermé numériquement effectif relativement à $\Delta$. Dans un premier temps, nous montrerons que ce courant est porté par le diviseur diagonal, à l'aide d'un lemme de feuilletage de son support similaire à celui établi par J. Duval dans \cite{Green presque complexe}. Dans un second temps, nous obtiendrons une contradiction aux points doubles par un argument d'homologie inspiré de M.L. McQuillan \cite{Bloch hyperbolicity}. Les deux étapes de la preuve feront abondamment appel à la théorie de recouvrement des surfaces d'Ahlfors pour les applications quasiconformes.\\

\noindent {\bf Remerciements:} Cet article doit énormément aux idées et aux nombreux conseils de Julien Duval. Qu'il en soit chaleureusement remercié.

\section{Préliminaires}

Ce chapitre énumère des définitions et des résultats généraux à propos de la géométrie du plan projectif presque complexe, des applications quasiconformes et plus particulièrement de la théorie de recouvrement des surfaces d'Ahlfors, ainsi que des courants positifs et des fonctions plurisousharmoniques dans un contexte presque complexe. Il s'agit principalement de rappels.

\subsection{Géométrie du plan projectif presque complexe}

La positivité des intersections de deux $J$-disques permet d'énoncer un résultat de type Hurwitz sur les limites de $J$-disques. Le second paragraphe est consacré au redressement de la structure presque complexe le long d'une $J$-droite, et nous sera utile au patagraphe 3.2 pour établir un lemme de Brunella presque complexe. Enfin, la notion d'éclaté presque complexe est rappelée en prévision chapitre 4.

\subsubsection{Positivité d'intersection}

Dans toute variété munie d'une structure presque complexe, les éventuelles intersections de deux $J$-disques d'images distinctes sont isolées et strictement positives en homologie (voir \cite{Gromov} et l'article de McDuff dans \cite{Audin}). Ceci entraîne qu'un $J$-disque $f:\mathbb{D}\rightarrow \mathbb{P}^{2}(\Cplx )$ limite d'une suite $(f_{n})$ de $J$-disques évitant une $J$-droite $L$ évite $L$ ou y est contenu.

\subsubsection{Redressement de la structure presque complexe le long d'une $J$-droite}

Soit $L$ une $J$-droite fixée. Il existe un difféomorphisme $\Phi$, défini près de $L$, et envoyant $L$ sur une droite complexe (une $i$-droite) le long de laquelle $\Phi_{*}J=d\Phi\circ J \circ d\Phi^{-1}$ coïncide avec $i$. En effet, d'après J.C. Sikorav \cite{Sikorav}, il existe une difféomorphisme de $\mathbb{P}^{2}(\Cplx )$ envoyant $L$ sur une droite complexe (une $i$-droite), et redressant tangentiellement $J$ en la structure standard $i$. Il est donc toujours possible de supposer que $L$ est une droite complexe, et que $J=i$ sur $TL$. Pour faire coïncider $J$ avec $i$ normalement le long de $L$, fixons un pinceau de $J$-droites $\mathcal{L}_{P}$ centré en un point $P$ n'appartenant pas à $L$; il existe un voisinage $U$ de $L$ et un difféomorphisme de $U$, fixant $L$, et envoyant tout $J$-disque $L'\cap U$ ($L'\in \mathcal{L}_{P}$) sur un disque holomorphe (un $i$-disque). Il est alors aisé de construire le difféomorphisme $\Phi$.\\
Ceci nous sera utile dans la démonstration de la proposition \ref{Brunella}.

\subsubsection{Eclatement presque complexe}

Rappelons la définition de l'éclaté du plan projectif presque complexe en un point (voir par exemple \cite{Green presque complexe}).\\
Par tout point $P$ du plan projectif presque complexe passe un pinceau de $J$-droites paramétré par $\mathbb{P}^{1}(\Cplx )$, donnant une projection centrale $\pi_{P}:\mathbb{P}^{2}(\Cplx )\setminus \{ P \} \rightarrow \mathbb{P}^{1}(\Cplx )$ (voir \cite{Gromov} et \cite{Sikorav}).\\
Ce pinceau se redresse localement sur le pinceau des droites complexes de $\mathbb{C}^{2}$ en $0$ à l'aide d'un difféomorphisme $\Psi_{P}$ défini près de $P$, transportant la structure $J$ en une structure presque complexe, encore notée $J$, coïncidant avec $i$ en $0$. Ce difféomorphisme est de classe $\mathcal{C}^{\infty}$ hors de $P$ mais seulement $\mathcal{C}^{1+Lip}$ en $P$. L'éclaté $\widetilde{X}_{P}$ de $(\mathbb{P}^{2}(\Cplx ),J)$ en $P$ est défini comme l'éclaté complexe usuel $\widetilde{\mathbb{C}^{2}}$ en $0$ via ce redressement local.\\

La projection $\pi_{P}$ se relève en une fibration $\widetilde{\pi}_{P}: \widetilde{X}_{P}\rightarrow \mathbb{P}^{1}(\Cplx )$. La structure presque complexe $\widetilde{J}_{P}$, définie hors de $E_{P}$ par relèvement de $J$, admet un prolongement lipschitz à $\widetilde{X}_{P}$ (voir \cite{Green presque complexe}), coïncidant avec $i$ sur $E_{P}$. Soit alors une $2$-forme $\lambda$ à support dans un voisinage de $E_{P}$ et positive pour $i$ tangentiellement à $E_{P}$: si $c>0$ est suffisament petit, la $2$-forme  $\widetilde{\omega}_{P}=\widetilde{\pi}_{P}^{*}\omega + c\lambda$ est non dégénérée sur $\widetilde{X}_{P}$ et positive pour $\widetilde{J}_{P}$, ie $\widetilde{\omega}_{P}(.,\widetilde{J}_{P}.)>0$.\\

A cause de la perte de régularité due au caractère non lisse de $\Psi_{P}$ en $P$, il n'est plus possible de définir des tours d'éclatements comme en complexe. Cependant, l'éclaté du plan projectif presque complexe en un nombre fini de points distincts est parfaitement défini.\\
Nous noterons $\widetilde{X}$ l'éclaté de $(\mathbb{P}^{2}(\Cplx ) , J)$ en les six points doubles $S_{k}$, $1\leq k \leq 6$, $\Pi :\widetilde{X}\rightarrow \mathbb{P}^{2}(\Cplx )$ la projection canonique contractant les diviseurs $E_{S_{k}}$ sur les points $S_{k}$, ainsi que $\widetilde{J}$ la structure presque complexe obtenue par prolongement de $J$ à $\widetilde{X}$. Nous fixerons une $2$-forme $\widetilde{\omega}$ non dégénérée et positive par rapport $\widetilde{J}$. \\

Etant donnée une $J$-droite $L$, l'adhérence $\widetilde{L}$ de la $\widetilde{J}$-courbe $\Pi^{-1}(L\setminus \{ S_{k}, 1\leq k \leq 6 \})$  est encore une $\widetilde{J}$-courbe de $\widetilde{X}$ appelée transformée stricte de $L$. Enfin, un $J$-disque $f$ évitant les points doubles se relève en un $\widetilde{J}$-disque $\widetilde{f}=\Pi^{-1}\circ f$ évitant les diviseurs exceptionnels.\\

L'homologie de $\widetilde{X}$ étant identique à celle de l'éclaté standard du plan projectif complexe en les points doubles, les transformées strictes $\widetilde{\Delta}_{1}$, $\widetilde{\Delta}_{2}$, $\widetilde{\Delta}_{3}$ des trois diagonales sont, comme les diviseurs exceptionnels, d'autointersection $-1$.\\

\subsection{Projections centrales et applications quasiconformes}

La projection centrale d'un $J$-disque induit une application du disque unité vers $\mathbb{P}^{1}(\Cplx )$ qui n'est pas en général holomorphe, mais quasiconforme, comme nous le rappelons dans le paragraphe suivant. Le lemme énoncé au paragraphe 2.2.2 sera utile à deux reprises, aux chapitres 3 et 4. Enfin, la théorie de reouvrement des surfaces d'Ahlfors est au coeur de la démonstration du théorème de Bloch. 

\subsubsection{Applications quasiconformes}
Nous renvoyons à \cite{Lehto} pour les propriétés énoncées dans ce paragraphe. Commençons par rappeler la définition d'une application quasiconforme. Soit $0< \alpha <1$. Une application continue $\varphi : \mathbb{D}\rightarrow \mathbb{P}^{1}(\Cplx )$ est dite $\alpha$-quasiconforme si $||\overline{\partial}\varphi ||\leq \alpha ||\partial \varphi ||$ presque partout, $\partial \varphi$ et $\overline{\partial}\varphi$ désignant respectivement les composantes $\Cplx$-linéaire et $\Cplx$-antilinéaire de la dérivée de $\varphi$ au sens des distributions. Le théorème d'Ahlfors-Bers permet d'écrire toute application $\alpha$-quasiconforme comme la composée $\phi \circ h$ d'une fonction holomorphe $\phi:\mathbb{D}\rightarrow \mathbb{P}^{1}(\Cplx )$ et d'un homéomorphisme $\alpha$-quasiconforme du disque unité. Les propriétés topologiques des fonctions holomorphes (et par suite la théorie de recouvrement des surfaces d'Ahlfors, voir paragraphe suivant) s'étendent donc aux applications quasiconformes. En outre, une suite d'homéomorphismes $\alpha$-quasiconformes du disque unité admet une sous-suite convergeant vers un homéomorphisme $\alpha$-quasiconforme, donc une suite d'applications $\alpha$-quasiconformes à valeurs dans un compact de $\mathbb{C}$ converge uniformément, après extraction d'une sous-suite.\\

Comme annoncé plus haut, d'après \cite{Green presque complexe}, pour tout point $P$ du plan projectif presque complexe, il existe une constante $0<\alpha_{P} <1$ telle que pour tout $J$-disque $f:\mathbb{D}\rightarrow \mathbb{P}^{2}(\Cplx )$ évitant $P$, l'application $\pi_{P}\circ f$ soit $\alpha_{P}$-quasiconforme. Si $f$ évite les six points doubles $S_{k}$ de la configuration $C$, les applications $v\pi_{S_{k}}\circ f$ sont donc $\alpha$-quasiconformes, avec $\alpha=\displaystyle{\max_{1\leq k \leq 6}}\alpha_{S_{k}}$.

\subsubsection{Critère de normalité pour une suite de $J$-disques évitant $C$}

La proposition suivante fournit un critère de normalité pour une suite de $J$-disques évitant la configuration $C$ et ne convergeant pas vers les diagonales, et découle directement des propriétés de convergence des suites d'applications quasiconformes.
\begin{Prop}
Soit $(f_{n})$ une suite de $J$-disques évitant $C$ et telle que $(f_{n}(\mathbb{D}))$ ne converge pas au sens de Hausdorff vers les diagonales de $C$. Alors si pour un $k\in[1,6]$ la suite d'applications $\alpha$-quasiconforme $\varphi_{k,n}=\pi_{S_{k}}\circ f_{n}$ admet une sous-suite convergente, $(f_{n})$ est normale.
\end{Prop}

\begin{demo}
Pour fixer les idées, prenons $k=1$, $S_{1}$ étant le point d'intersection de $L_{1}$ et $L_{2}$. Supposons que $(\varphi_{1,n})$ converge, après extraction d'une sous-suite, vers une application $\alpha$-quasiconforme $\varphi$. Une $J$-courbe entière $f:\Cplx \rightarrow \mathbb{P}^{2}(\Cplx )$ obtenue après reparamétrisation à la Brody de $f_{n}$ au voisinage d'un éventuel point d'explosion $e$ est contenue dans la $J$-droite $L$ du pinceau centré en $S_{1}$ située au-dessus de $\varphi (e)$. D'une part, par positivité d'intersection, $f(\Cplx )$ évite $C$, et d'autre part, comme $f$ est holomorphe pour la structure complexe induite sur $L$ par $J$, le petit théorème de Picard interdit que $f$ évite plus de deux points; il en découle que $L$ est la diagonale de la configuration $C$ passant par $S_{1}$. Comme $\varphi$ ne peut être constante égale à $\varphi (e)$, ce qui impliquerait que $(f_{n}(\mathbb{D}))$ converge au sens de Hausdorff vers un diagonale, l'ensemble $E$ des points d'explosion est discret (principe des zéros isolés pour une application quasiconforme).\\
Observons la suite $(f_{n})$ depuis un deuxième point double, par exemple le point d'intersection $S_{3}$ de $L_{1}$ et $L_{3}$. La suite d'applications $\alpha$-quasiconformes $(\varphi_{3,n})$ converge uniformément, après extraction d'une sous-suite, sur tout compact de $\mathbb{D}\setminus E$. D'après le principe du maximum, elle converge uniformément sur tout compact de $\mathbb{D}$. Comme les projections $\pi_{S_{1}}$ et $\pi_{S_{3}}$ fournissent des coordonnées sur $\mathbb{P}^{2}(\Cplx )\setminus L_{1}$, l'ensemble $E$ est vide et $(f_{n})$ est normale. \qed
\end{demo}

\subsubsection{Théorie de recouvrement des surfaces de Riemann}

Ce paragraphe rappelle les résultats essentiels de la théorie de recouvrement des surfaces de Riemann due à Ahlfors, qui permet de généraliser les principaux énoncés de la théorie de distribution des valeurs des fonctions holomorphes aux applications quasiconformes.\\
Fixons une surface de Riemann $\Sigma_{0}$ munie d'une métrique hermitienne $g$. Soit $\Sigma$ une surface de Riemann compacte connexe à bord lisse, et soit $\varphi :\overline{\Sigma} \rightarrow \Sigma_{0}$ une fonction lisse sur $\Sigma$ et $\alpha$-quasiconforme sur l'intérieur de $\Sigma$. Nous noterons $v_{\varphi^{*}g}$ la forme d'aire de $\varphi^{*}g$.  Posons $\operatorname{Aire}(\varphi (\Sigma ))=\int_{\Sigma}v_{\varphi ^{*}g}$ l'aire avec multiplicité de $\varphi (\Sigma )$ et $L=\int_{\partial \Sigma \setminus \varphi^{-1}(\partial \Sigma_{0})}||d\varphi ||_{g}$ la longueur de son bord relatif. Si $\Omega$ est un domaine compact à bord lisse de $\Sigma_{0}$, notons $F(\Omega )=\frac{\operatorname{Aire}(\varphi (\varphi^{-1}(\Omega)))}{\operatorname{Aire}(\Omega )}$ le nombre moyen de feuillets au-dessus de $\Omega$.\\
Le théorème suivant peut s'interpréter comme une formulation géométrique du premier théorème principal de Nevanlinna; il affirme qu'à un terme de bord près, le nombre moyen de feuillets au-dessus d'un domaine est égal à $F(\Sigma_{0})$.

\begin{Theo}\label{recouvrement}
Soit $\Omega$ un domaine compact à bord lisse de $\Sigma_{0}$. Il existe une constante $h>0$ ne dépendant ni de $\Sigma$ ni de $\varphi $ telle que $|F(\Omega ) - F(\Sigma_{0})|\leq hL$.
\end{Theo}
Par la suite, nous appellerons $h$ toute constante indépendante de $\varphi $.
L'inégalité suivante, dite inégalité d'Ahlfors, constitue le résultat central de la théorie de recouvrement des surfaces de Riemann:

\begin{Theo}\label{inégalité}
Il existe une constante positive $h$ ne dépendant que de $(\Sigma_{0}, g)$ telle que:
$$\min(0, \chi (\Sigma ))\leq F(\Sigma_{0})\chi (\Sigma_{0}) + hL$$
où les entiers $\chi (\Sigma )$ et $\chi (\Sigma_{0})$ désignent les caractéristiques d'Euler-Poincaré des surfaces $\Sigma$ et $\Sigma_{0}$.
\end{Theo}

Ces énoncés permettent d'obtenir un second théorème principal. Munissons $\mathbb{P}^{1}(\Cplx )$ de la métrique associée à la forme de Fubini-Study normaliséee $\omega '$. Soient $\Omega_{1}$, $\Omega_{2}$ et $\Omega_{3}$ des disques ouverts de $\mathbb{P}^{1}(\Cplx )$ d'adhérences deux à deux disjointes. Considérons une composante connexe $\Omega$ de $\varphi^{-1}(\Omega_{i})$, $1\leq i \leq 3$. Si $\Omega$ est relativement compacte dans $\Sigma$, la restriction de $\varphi $ à $\Omega$ possède un degré. On dit alors que $\Omega$ est un île au-dessus de $\Omega_{i}$. Dans le cas contraire, $\Omega$ est appelée péninsule au-dessus de $\Omega_{i}$.\\
Posons alors $\Sigma_{0}=\mathbb{P}^{1}(\Cplx )\setminus (\Omega_{1} \cup \Omega_{2} \cup \Omega_{3})$ et $\Sigma '= \varphi^{-1}(\Sigma_{0})$. Si $I$ désigne le nombre d'îles au-dessus des disques $\Omega_{1}$, $\Omega_{2}$ et $\Omega_{3}$, alors $\chi (\Sigma )- \chi (\Sigma ') = I$. \\
Par le théorème \ref{inégalité}, il existe une constante $h>$ telle que $F(\Sigma_{0}) \leq -\min (0,\chi (\Sigma ) -I) + hL$. Mais par ailleurs, d'après le théorème \ref{recouvrement}, $\operatorname{Aire}(\varphi (\Sigma ))\leq F(\Sigma_{0}) + hL$ ce qui permet d'obtenir l'inégalité:
\begin{equation}\label{eq1}
\operatorname{Aire}(\varphi (\Sigma )) \leq -\min (0,\chi (\Sigma )-I) + hL \leq I - \min (0,\chi (\Sigma )) + hL.
\end{equation}

Nous nous intéresserons particulièrement au cas où $\varphi$ évite les valeurs $0$ et $\infty$. Supposons alors que $0\in \Omega_{1}$ et $\infty \in \Omega_{2}$. Il ne peut y avoir d'îles au-dessus de $\Omega_{1}$ et $\Omega_{2}$: $I$ correspond alors au nombre d'îles au-dessus de $\Omega_{3}$.\\

Nous pouvons encore préciser l'inégalité (\ref{eq1}) en vue du paragraphe $3$. Le domaine $\Sigma$ est alors un disque, de caractéristique d'Euler-Poincaré égale à $1$. Notons $I_{1}$ le nombre d'îles de degré $1$ au-dessus de $\Omega_{3}$. Le théorème \ref{recouvrement} ainsi que (\ref{eq1}) permettent d'obtenir la double inégalité $I_{1} + 2(I-I_{1}) \leq F(\Omega_{3}) \leq I_{1} + (I-I_{1}) + hL$, soit $I\leq I_{1} + hL$. On dispose donc d'un second théorème principal dont le second membre ne prend en compte que les îles de degré $1$ au-dessus de $\Omega_{3}$:
\begin{equation}\label{eq2}
\operatorname{Aire}(\varphi (\Sigma )) \leq I_{1} + hL.
\end{equation}
 
Ces inégalités (\ref{eq1}) et (\ref{eq2}) seront utilisées sous forme intégrée "à la Nevanlinna"; les termes de longueur et de caractéristique d'Euler-Poincaré seront négligeables devant le terme d'aire.

\subsection{Courants positifs et fonctions plurisousharmoniques}

Nous renvoyons à \cite{Pali} pour les résultats énoncés dans ce paragraphe, excepté pour le théorème de restriction des courants positifs fermés à une $J$-droite, que l'on trouvera dans \cite{Elkhadhra}.

\subsubsection{Courants positifs}

L'automorphisme $J$ se prolonge en un automorphisme $\Cplx$-linéaire, encore noté $J$, du complexifié $\Cplx \otimes_{\mathbb{R}}T_{\mathbb{R}}\mathbb{P}^{2}(\Cplx )$ du fibré tangent réel $T_{\mathbb{R}}\mathbb{P}^{2}(\Cplx )$ du plan projectif. Notons respectivement $T_{J}^{1,0}$ et $T_{J}^{0,1}$ les sous-fibrés de $\Cplx \otimes_{\mathbb{R}}T_{\mathbb{R}}\mathbb{P}^{2}(\Cplx )$ dont les fibres sont les sous-espaces propres de $J$ associés aux valeurs propres $i$ et $-i$,  et $T_{J,1,0}^{*}$,$T_{J,0,1}^{*}$ leurs duaux. Notons enfin $\bigwedge_{J}^{1,1}$ l'espace des formes différentielles de classe $\mathcal{C}^{\infty}$, de bidegré $(1,1)$.\\
Un courant $T$ est dit de type $(1,1)$ s'il s'écrit localement comme une forme différentielle de bidegré $(1,1)$ dont les coefficients sont des distributions.\\

Une forme $\theta$ de bidegré $(1,1)$ est dite positive si pour tout champ de vecteur $\zeta$, $\theta (\zeta,J\zeta)\geq 0$. Comme nous sommes en bidegré $(1,1)$, une telle forme différentielle s'écrit localement comme une somme de formes du type $i\alpha \wedge \overline{\alpha}$, où $\alpha \in T_{J,1,0}^{*}$.\\
Un courant $T$ de bidegré $(1,1)$ est positif si pour tout $(1,1)$-forme $\theta$ positive, $T(\theta )\geq 0$. Si $(\zeta_{1}, \zeta_{2})$ est un repère local du fibré $T_{J}^{1,0}$, d'après \cite{Pali}, $T$ s'écrit dans le repère dual $(\zeta_{1}^{*},\zeta_{2}^{*})$:
$$T=iT_{1,1}\zeta_{1}^{*}\wedge \overline{\zeta}_{1}^{*} + iT_{1,2}\zeta_{1}^{*}\wedge \overline{\zeta}_{2}^{*} + iT_{2,1}\zeta_{2}^{*}\wedge \overline{\zeta}_{1}^{*} + iT_{2,2}\zeta_{2}^{*}\wedge \overline{\zeta}_{2}^{*}$$
où les $T_{r,s}$ sont des distributions d'orde zéro, donc des mesures de Radon complexes, avec $T_{1,1}$ et $T_{2,2}$ positives. De plus, $||T_{1,2}||, ||T_{2,1}|| \leq 2T_{1,1} + 2T_{2,2}$. En particulier, en écrivant la forme locale du courant dans le repère $(\lambda_{1}\zeta_{1}, \lambda_{2}\zeta_{2})$, avec $\lambda_{2}\gg \lambda_{1}>0$, on on montre que si $T_{2,2}=0$, alors $T_{1,2}=T_{2,1}=0$.\\

Les courants positifs sont d'ordre zéro, et une suite de courants positifs de masse uniformément bornée admet donc une valeur d'adhérence pour la topologie faible, qui est encore un courant positif. Si $c$ est une courbe $J$-holomorphe, le courant $[c]$ d'intégration sur $c$ est positif.\\

Bien que ne disposant pas d'un théorème de Siu pour les courants positifs fermés, nous pouvons malgré tout isoler leur composante singulière portée par une $J$-droite donnée: si $T$ est un courant positif fermé du plan projectif presque complexe, pour toute $J$-droite $L$, le courant $\mathds{1}_{L}T$ est positif fermé (voir \cite{Elkhadhra}), et s'écrit donc  $\mathds{1}_{L}T=\alpha [L]$, où $\alpha$ est un réel positif. En particulier, $T$ admet la décomposition suivante:
$$T = T' + \alpha_{1}[\Delta_{1}] + \alpha_{2}[\Delta_{2}] + \alpha_{3}[\Delta_{3}]$$ 
où les $\alpha_{k}$ sont positifs et $T'$ n'admet aucune diagonale comme composante singulière. 

\subsubsection{Fonctions plurisousharmoniques et voisinages pseudoconvexes} 

Rappelons brièvement ce qu'est une fonction plurisousharmonique dans le contexte presque complexe.\\
Si $U$ est un ouvert de $\mathbb{P}^{2}(\Cplx )$, une application $u:U\rightarrow \mathbb{R}\cup \{ -\infty \}$ est dite plurisousharmonique si pour tout $J$-disque $f:\mathbb{D}\rightarrow U$, l'application $u\circ f$ est sous-harmonique. \\

Si $P$ est un point du plan projectif presque complexe, il existe une difféomorphisme $u$ d'un voisinage de $P$ dans un ouvert de $\mathbb{C}^{2}$ envoyant $J_{P}$ sur $i$. D'après \cite{Pali}, $P\mapsto ||u(P)||^{2}$ ($||.||$ désignant ici la norme hermitienne de $\Cplx^{2}$) est plurisousharmonique. Tout point $P$ possède donc une base $(U_{t})_{0\leq t\leq \varepsilon}$ de voisinages pseudo-convexes, donnés par: $U_{t}=u^{-1}([0,t])$. Ainsi, le diviseur exceptionnel $E_{P}$ au-dessus de $P$ dans l'éclaté presque complexe de $\widetilde{X}_{P}$ de $V$ en $P$ possède lui aussi une base de voisinages pseudo-convexes $(\widetilde{U}_{t}=\widetilde{\pi}_{P}^{-1}(U_{t}))_{0\leq t\leq \varepsilon}$. \\

Si $f:\overline{\mathbb{D}}\rightarrow \mathbb{P}^{2}(\Cplx )$ est un $J$-disque évitant $P$ et dont le relevé $\widetilde{f}_{P}$ dans $\widetilde{X}_{P}$ vérifie: $\widetilde{f}_{P}(\partial \mathbb{D})\subset \partial \widetilde{U}_{t}$ pour un certain $0<t<\varepsilon$, le principe du maximum pour les fonctions sous-harmoniques interdit que $\widetilde{f}_{P}(\mathbb{D})$ soit contenu dans $\widetilde{U}_{t'}\setminus \widetilde{U}_{t}$, pour un certain $t<t'<\varepsilon$. Cette remarque nous sera utile au dernier paragraphe.\\

Venons-en à la démonstration proprement dite du théorème de Bloch.

\section{Construction de courants de Nevanlinna}

Donnons-nous une suite $(f_{n})$ non normale de $J$-disques. Raisonnons par l'absurde, et supposons que la suite $f_{n}(\mathbb{D})$ ne converge pas vers $\Delta$ au sens de Hausdorff. Il est toujours possible de supposer que $f_{n}(0)$ reste loin de $\Delta$:  quitte à reparamétrer les J-disques par un automorphisme du disque unité et à extraire une sous-suite de $(f_{n})$, nous pouvons supposer que pour un certain $\varepsilon_{0} >0$, $d(f_{n}(0), \Delta_{0} ) \geq \varepsilon_{0}$. Notons que malgré cette supposition, les courbes entières issues des explosions de la suite des J-disques peuvent quand-même être portées par les diagonales (voir la remarque sur le théorème de Brody dans l'introduction).\\
Ce paragraphe est dédié à la construction,  à partir des suites $(f_{n})$ et $(\widetilde{f}_{n})$, de courants $T$ et $\widetilde{T}$ de bidegré $(1,1)$, positifs, fermés et numériquement effectifs (nef). Classiquement, nous allons les construire à partir des fonctions caractéristiques de Nevanlinna (paragraphe 3.1). L'effectivité numérique se déduit d'un Lemme de Brunella presque complexe (paragraphe 3.2), ce qui justifie l'utilisation de courants de Nevanlinna et non de courants d'Ahlfors.

\subsection{Construction de courants fermés}

Rappelons les définitions des fonctions caractéristiques de Nevanlinna. Soit $(M,J,\Theta )$ une variété presque complexe de dimension réelle $4$ munie d'une $(1,1)$-forme positive $\Theta$ induisant une métrique hermitienne $||.||_{\Theta}$ (Dans le contexte, $(M,J, \Theta )=(\mathbb{P}^{2}(\Cplx ),J, \omega )$ ou $(\widetilde{X}, \widetilde{J}, \widetilde{\omega})$. Si $R$ est une réel strictement positif et $g:D(0,R)\rightarrow M$ une application J-holomorphe, notons pour $0<r<R$:
$$T_{g,r}=\int_{0}^{r}[g(D(0,t))]\frac{dt}{t} \hspace{2cm} \text{ le courant caractéristique de } g$$
$$T_{g,r}(\Theta )=\int_{0}^{r}\left (\int_{D(0,t)}g^{*}\Theta \right ) \frac{dt}{t} \hspace{0.6cm} \text{ la caractéristique d'aire de } g$$ 
$$L_{g,r}=\int_{0}^{r}\int_{\partial D(0,t)}||dg||_{\Theta}\frac{dt}{t} \hspace{2.1cm} \text{ la caractéristique de longueur de } g.$$

Le courant $T_{g,r}$ est positif (voir 2.3.1). Nous recherchons une suite $(r_{n})$ telle que: $\displaystyle{\lim_{n\rightarrow +\infty}}\frac{L_{f_{n},r_{n}}}{T_{f_{n},r_{n}}(\omega )}=\displaystyle{\lim_{n\rightarrow +\infty}}\frac{L_{\widetilde{f}_{n},r_{n}}}{T_{\widetilde{f}_{n},r_{n}}(\widetilde{\omega} )}$; les suites de courants positifs $\left( \frac{T_{f_{n},r_{n}}}{T_{f_{n},r_{n}(\omega )}} \right)$ et $\left( \frac{T_{\widetilde{f}_{n},r_{n}}}{T_{\widetilde{f}_{n},r_{n}(\widetilde{\omega} )}} \right)$ convergeront alors, après extraction d'une sous-suite, vers des courants fermés. \\
Les caractéristique d'aire $T_{\varphi ,r}(\omega ')=\int_{0}^{r}\left (\int_{D(0,t)}\varphi^{*}\omega ' \right ) \frac{dt}{t}$  et de longueur $L_{\varphi,r}=\int_{0}^{r}\int_{\partial D(0,t)}||d\varphi ||_{\omega '}\frac{dt}{t}$ d'une application quasiconforme $\varphi :D(0,R)\rightarrow \mathbb{P}^{1}(\Cplx )$ relativement à la forme de Fubini-Study normalisée $\omega '$ sont bien définies; en vue du paragraphe 3, la suite $(r_{n})$ recherchée devra en outre vérifier: $\displaystyle{\lim_{n\rightarrow +\infty}}\frac{L_{\varphi_{k,n},r_{n}}}{T_{\varphi_{k,n},r_{n}}(\omega ')}=0$.\\

Classiquement, l'existence d'une telle suite $(r_{n})$ repose sur l'inégalité longueur-aire suivante (voir par exemple \cite{Brunella}):
\begin{equation}\label{longueur-aire}
\forall 0<\eta <r<R, \hspace{1cm} \frac{L_{g,r}-L_{g,\eta}}{T_{g,r}(\Theta )}\leq \sqrt{2\pi r \log (\frac{r}{\eta})\frac{1}{T_{g,r}(\Theta )^{2}}\frac{dT_{g,r}(\Theta )}{dr}}.
\end{equation}
Pour une courbe entière $f:\Cplx \rightarrow \mathbb{P}^{2}(\Cplx )$ non constante, cette inégalité permet immédiatement de trouver une suite $(R_{n})$ croissant vers l'infini telle que $\displaystyle{\lim_{n\rightarrow +\infty}}\frac{L_{f,R_{n}}}{T_{f,R_{n}}(\omega )}=0$. En effet, la fonction positive $\frac{2\pi r \log (r)}{T_{f,r}(\omega )^{2}}\frac{dT_{f,r}(\omega )}{dr}$ est intégrable sur $[2, +\infty [$ par rapport à la mesure de masse infinie $\frac{dr}{r\log (r)}$:
$$\int_{2}^{R}\frac{2\pi r\log (r)}{T_{f,r}(\omega )^{2}}\frac{dT_{f,r}(\omega )}{dr} \frac{dr}{2\pi \log (r)}=\frac{1}{T_{f,2}(\omega )}- \frac{1}{T_{f,R}(\omega )} < +\infty .$$
Pour tout $\delta >0$, l'ensemble des $r\geq 2$ tels que $\frac{L_{f,r}-L_{f,2}}{T_{f,r}(\omega )}>\delta $ est donc de mesure de Lebesgue finie. Comme par ailleurs $\displaystyle{\lim_{r\rightarrow  +\infty}}T_{f,r}(\omega )=+\infty$, la suite $(R_{n})$ est aisément construite.\\

Le cas de la suite de J-disques $(f_{n}(\mathbb{D}))$ est moins immédiat. Remarquons d'abord que la suite $\operatorname{Aire}(f_{n}(D(0,r)))$ est non bornée pour un certain $0<r<1$. Il suffit pour cela de voir que les explosions de $(f_{n})$ sont d'aire infinie: soit $f:\Cplx\rightarrow \mathbb{P}^{2}(\Cplx )$ une $J$-courbe de Brody obtenue par reparamétrisation en un point d'explosion $e$ de la suite $(f_{n})$. Si $f(\Cplx )$ était d'aire finie, $f$ se prolongerait en une $J$-courbe rationnelle $\overline{f}:\mathbb{P}^{1}(\Cplx )\rightarrow \mathbb{P}^{2}(\Cplx )$ appelée bulle (voir l'article de J.-C. Sikorav dans \cite{Audin}). Comme l'intersection d'une $J$-droite et d'une $J$-courbe rationnelle non constante n'est jamais vide (voir par exemple \cite{Green presque complexe} pour la démonstration d'un théorème de Liouville presque complexe), cette bulle rencontrerait $C$ en au moins deux points distincts. La $J$-courbe entière $f(\Cplx )=\overline{f}(\mathbb{P}^{1}(\Cplx )\setminus \{\infty \})$, couperait donc $C$ en au moins un point, et par positivité d'intersection, serait contenue dans l'une des quatre $J$-droites $L_{j}$ tout en évitant les trois autres, contredisant le petit théorème de Picard. \\

Ainsi, pour un certain $0<r_{0}<1$, $\displaystyle{\limsup_{n\rightarrow +\infty}T_{f_{n},r_{0}}}(\omega ) =+\infty$. Si $\eta$ est compris strictement entre $0$ et $r_{0}$, pour tout $r>r_{0}$ on a:
$$\lim_{n\rightarrow +\infty}\int_{r_{0}}^{r}\frac{2\pi t\log (t/\eta )}{T_{f_{n},t}(\omega )^{2}}\frac{dT_{f_{n},t}(\omega )}{dt} \frac{dt}{2\pi \log (t/\eta )}= \lim_{n\rightarrow +\infty} \left( \frac{1}{T_{f_{n},r_{0}}(\omega )}-\frac{1}{T_{f_{n},r}(\omega )} \right) = 0$$
donc pour presque tout $r>r_{0}$, $\displaystyle{\liminf_{n\rightarrow +\infty}}\frac{L_{f_{n},r}-L_{f_{n},\eta}}{T_{f_{n},r}(\omega )}=0$. Ce n'est pas suffisant pour conclure, car rien n'indique a priori que $\displaystyle{\liminf_{n\rightarrow +\infty}}\frac{L_{f_{n},\eta}}{T_{f_{n},r}(\omega )}=0$. Notons alors $D(0,\eta_{n})$ le disque centré en $0$ dont l'image par $f_{n}$ est d'aire $1$. La suite des J-disques $f_{n}(\overline{D(0,\eta_{n} )})$ est normale, puisqu'elle est d'aire bornée et évite $C$, donc la suite $(L_{f_{n},\eta_{n}})$ est bornée. Pour tout $r>r_{0}$, $\frac{L_{f_{n},\eta_{n}}}{T_{f_{n},r}(\omega )}$ a donc pour limite inférieure $0$. Distinguons deux cas:

\begin{itemize}
\item Si $(\eta_{n})$ est minorée par un réel strictement positif $\eta$, alors pour presque tout $r>r_{0}$, $\frac{L_{f_{n},r }}{T_{f_{n},r}(\omega )}$ converge vers $0$, et la suite constante $r_{n}=r$ convient.\\
\item Si la suite $(\eta_{n})$ tend vers $0$, il est tentant de reparamétrer l'anneau de grand module $\overline{D(0,r)}\setminus D(0,\eta_{n})$ par l'anneau $\overline{D(0,r/\eta_{n})}\setminus D(0,1)$, de façon à nous ramener à une situation similaire à celle d'une courbe entière. Cela équivaut à effectuer le changement de variable suivant:
$$\int_{2\eta_{n}}^{r}\frac{\varphi_{n}(t)dt}{t\log (t/\eta_{n})}=\int_{2}^{r/\eta_{n}}\frac{\varphi_{n}(\eta_{n}u)du}{u\log (u)}\hspace{0.5cm} \text{ avec } \varphi_{n}(t)=\frac{t\log (t/\eta_{n})}{T_{f_{n},t}(\omega )}\frac{dT_{f_{n},t}(\omega )}{dt}.$$
Cette intégrale vaut $1/T_{f_{n},2\eta_{n}}(\omega ) - 1/T_{f_{n},r}(\omega )$, et est donc majorée par $\frac{1}{T_{f_{n},\eta_{n}}(\omega )}=1$. Comme la fonction $\frac{1}{u\log u}$ n'est pas intégrable sur $[2,+\infty[$, pour tout $\delta >0$, l'ensemble $\{  u\in [2,r/\eta_{n}] \text{ | } \varphi_{n}(\eta_{n}u)\geq \delta \}$ est de mesure de Lebesgue majorée uniformément en $n$ par une constante $M_{\delta}$. Nous en déduisons que l'ensemble $\{ t\in ]2\eta_{n},r] \text{ | } \varphi_{n}(t)\geq \delta \}$ est de mesure de Lebesgue inférieure à $\eta_{n}M_{\delta}$. Pour presque tout $0<r<1$, la suite $\frac{L_{f_{n},r}}{T_{f_{n},r}(\omega )}$ a donc pour limite inférieure $0$. 
\end{itemize}

Dans tous les cas, pour presque tout $r\in ]r_{0},1]$, $\displaystyle{\liminf_{n\rightarrow +\infty}}\frac{L_{f_{n},r}}{T_{f_{n},r}(\omega )}=0$. Voyons les cas des suites $(\widetilde{f}_{n})$ et $(\varphi_{k,n})$. \\
La suite $(\operatorname{Aire}(\widetilde{f}_{n}(\mathbb{D}))$ n'est évidemment pas bornée, donc ce qui précède reste valable, et pour presque tout $r$ suffisamment proche de $1$: $\displaystyle{\liminf_{n\rightarrow +\infty}}\frac{L_{\widetilde{f}_{n},r}}{T_{\widetilde{f}_{n},r}(\widetilde{\omega})}=0$. \\
D'après le paragraphe 2.2.2, les suites $(\varphi_{k,n})$ ne sont pas normales non-plus. Comme elles évitent deux points de $\mathbb{P}^{1}(\Cplx )$, les suites $T_{\varphi_{k,n},r}(\omega ')$ ne sont pas bornées pour $r$ suffisament proche de 1 (l'argument de prolongement des courbes entières d'aire finie utilisé plus haut s'adapte sans problèmes aux applications quasiconformes). Le raisonnement précédent s'applique encore, et pour presque tout $r$ assez proche de 1: $\displaystyle{\liminf_{n\rightarrow +\infty}}\frac{L_{\varphi_{k,n},r}}{T_{\varphi_{k,n},r}(\omega_{1})}=0$.\\

Choisissons alors un rayon $r$ tel que $\displaystyle{\liminf_{n\rightarrow +\infty}}\frac{L_{f_{n},r}}{T_{f_{n},r}(\omega )}=\displaystyle{\liminf_{n\rightarrow +\infty}}\frac{L_{\widetilde{f}_{n},r}}{T_{\widetilde{f}_{n},r}(\widetilde{\omega} )}=\displaystyle{\liminf_{n\rightarrow +\infty}}\frac{L_{\varphi_{k,n},r}}{T_{\varphi_{k,n},r}(\omega ')}=0$ (pour $1\leq k\leq 6$). Les suites de courants $\frac{T_{f_{n},r}}{T_{f_{n},r}(\omega )}$ et $\frac{T_{\widetilde{f}_{n},r}}{T_{\widetilde{f}_{n},r}(\widetilde{\omega})}$ convergent faiblement (après éventuelle extraction d'une sous-suite) vers des courants fermés $T$ et $\widetilde{T}$. Quitte à reparamétrer le disque $D(0,r)$ par le disque unité, nous pouvons supposer que $r=1$. Les suites $T_{f_{n},r}$, $T_{\widetilde{f}_{n},r}$ et $T_{\varphi_{k,n},r}(\omega ')$ ainsi que $L_{f_{n},r}$, $L_{\widetilde{f}_{n},r}$ et $L_{\varphi_{k,n},r}$ seront notées $T_{f_{n}}$, $T_{\widetilde{f}_{n}}$ et $T_{\varphi_{k,n}}(\omega ')$ ainsi que $L_{f_{n}}$, $L_{\widetilde{f}_{n}}$ et $L_{\varphi_{k,n}}$.

\subsection{Intersection homologique et intersection géométrique}

L'objet de ce paragraphe est double: d'une part, rappeler que la caractéristique d'aire de $f_{n}$ domine asymptotiquement celle de $\varphi_{k,n}$ ($1\leq k \leq 6)$, ce qui peut s'interprêter comme un premier théorème principal de Nevanlinna, et d'autre part, minorer l'intersection homologique $[\widetilde{T}].[\widetilde{\Delta}_{l}]$ ($1\leq l \leq 3)$ par un terme d'intersection géométrique asymptotique. Les deux inégalités reposent sur la formule de Jensen, et s'obtiennent en mimant les techniques classiquement utilisées dans le contexte holomorphe.\\

Notons $\# [C].[C']$ le nombre de points d'intersections de deux $J$-courbes distinctes $C$ et $C'$, comptés avec multiplicité. Pour tout $n\geq 1$ et tout $1\leq l \leq 6$ posons:
$$I_{f_{n}}(\Delta_{l})=\int_{0}^{1}\# [f_{n}(D(0,t))].[\Delta_{l}]\frac{dt}{t} \hspace{0.5cm}\text{ et }\hspace{0.5cm} I_{\widetilde{f}_{n}}(\widetilde{\Delta}_{l})=\int_{0}^{1}\# [\widetilde{f}_{n}(D(0,t))].[\widetilde{\Delta}_{l}]\frac{dt}{t}.$$
Les intersections géométriques asymptotiques $I(T,\Delta_{l})$ et $I(\widetilde{T},\widetilde{\Delta}_{l})$ sont définies ainsi:
$$I(T,\Delta_{l})=\limsup_{n\rightarrow +\infty} \frac{I_{f_{n}}(\Delta_{l})}{T_{f_{n}}(\omega )} \hspace{0.5cm} \text{ et }\hspace{0.5cm} I(\widetilde{T},\widetilde{\Delta}_{l})=\limsup_{n\rightarrow +\infty} \frac{I_{\widetilde{f}_{n}}(\widetilde{\Delta}_{l})}{T_{\widetilde{f}_{n}}(\widetilde{\omega})}.$$

La première inégalité a déjà été établie et utilisée par Julien Duval dans \cite{Green presque complexe}. Elle nous sera utile au paragraphe suivant.
\begin{Prop}\label{comparaison projection}
Pour tout $1\leq k \leq 6$:
$$\limsup_{n\rightarrow +\infty} \frac{T_{\varphi_{k,n}}(\omega ')}{T_{f_{n}}(\omega )}\leq 1.$$
\end{Prop}

\begin{demo}
L'inégalité est une conséquence de l'existence d'une fonction $u$ négative à singularité logarithmique en $S_{k}$ et d'une $1$-forme bornée $\alpha$ définie sur $\mathbb{P}^{2}(\Cplx ) \setminus \{ S_{k} \}$ telles que:
\begin{equation}\label{logarithmique}
\pi_{S_{k}}^{*}\omega ' = \omega + dd_{J}^{c}u + d\alpha .
\end{equation}
En appliquant $T_{f_{n}}$ à cette égalité, nous obtenons, grâce au théorème de Stokes et à la formule de Jensen:
\begin{align}
T_{\varphi_{k,n}}(\omega ') &=T_{f_{n}}(\omega ) + \frac{1}{2\pi}\int_{0}^{2\pi}u\circ f_{n}(e^{i\theta})d\theta - u\circ f_{n}(0) + \int_{0}^{1}\int_{\partial \mathbb{D}}f_{n}^{*}\alpha \frac{dt}{t} \notag \\
&\ \leq T_{f_{n}}(\omega ) - u\circ f_{n}(0) + O(L_{f_{n}}). \notag
\end{align}
Comme $f_{n}(0)$ reste loin de $S_{k}$ pour tout $n$, $-u\circ f_{n}(0)$ est majoré. La proposition découle alors du caractère fermé de $T$.\\
Rappelons brièvement d'où provient l'identité (\ref{logarithmique}): les $2$-formes $\pi_{S_{k}}^{*}\omega '$ et $\omega$ étant cohomologues sur $\mathbb{P}^{2}(\Cplx )\setminus \{S_{k} \}$, il s'agit d'un problème local en $S_{k}$. Après redressement du pinceau de $J$-droites en $S_{k}$ par $\Psi_{S_{k}}$, $\pi_{S_{k}}^{*}\omega '$ s'écrit:
$\pi_{S_{k}}^{*}\omega ' = \omega + dd^{c}\log ||Z|| = \omega + dd_{J}^{c}\log ||Z|| + d(d\log ||Z||)\circ (J-i)$
où $Z$ est une coordonnée locale près de $S_{k}$ et $||.||$ la norme hermitienne de $\mathbb{C}^{2}$. Comme $\Psi_{S_{k}}$ est de classe $\mathcal{C}^{1+lip}$, $J(Z)-i=O(||Z||)$ et la $1$-forme $\alpha = d(\log ||Z||)\circ (J-i)$ est bornée. \qed
\end{demo}
Nous déduisons de cette proposition un permier théorème principal à la Nevanlinna: si $\Delta_{l}$ est la diagonale passant par $S_{k}$, le théorème \ref{recouvrement} appliqué à la suite d'applications $\varphi_{k,n}$ et au point $\pi_{S_{k}}(\Delta_{l})$ donne:
$$\limsup_{n\rightarrow + \infty}\frac{I_{f_{n}}(\Delta_{l})}{T_{\varphi_{k,n}(\omega ')}}\leq 1$$
ce qui entraîne:
$$I(T,\Delta_{l})\leq 1.$$
Cette dernière inégalité exprime en particulier le caractère numériquement effectif de $T$.\\

La seconde inégalité recherchée est un Lemme de Brunella presque complexe relatif au courant $\widetilde{T}$ et à $\widetilde{\Delta}_{l}$ ($1\leq l \leq 3$):
\begin{Prop}\label{Brunella}
$[\widetilde{T}].[\widetilde{\Delta}_{l}]\geq I(\widetilde{T}, \widetilde{\Delta}_{l})$.
\end{Prop}

\begin{demo}
Le procédé est similaire à celui utilisé pour démontrer la proposition \ref{comparaison projection}: nous allons nous ramener à une situation proche du cas holomorphe. L'inégalité proviendra alors du théorème de Lelong-Poincaré et de l'inégalité de Jensen, comme pour le lemme de Brunella (voir \cite{Brunella}).\\
Le problème est localisé au voisinage de $\widetilde{\Delta}_{l}$. D'après le paragraphe 2.1.2, nous pouvons supposer que $\Delta_{l}$ est une droite complexe le long de laquelle $J$ coïncide avec $i$. La structure relevée $\widetilde{J}$ coïncide alors avec $i$ le long de $\widetilde{\Delta}_{l}$.\\
Fixons une métrique hermitienne $|.|$ sur le fibré en droites associé au diviseur $\widetilde{\Delta}_{l}$, qui soit plate hors d'un petit voisinage de $\widetilde{\Delta}_{l}$. Si $s$ est une section de ce fibré s'annulant le long de $\widetilde{\Delta}_{l}$, d'après le théorème de Lelong-Poincaré, la $2$-forme $\omega_{\widetilde{\Delta}_{l}}$ définie au sens des courants par: 
\begin{equation}\label{Lelong-Poincaré}
\omega_{\widetilde{\Delta}_{l}}=[\widetilde{\Delta}_{l}] -dd^{c}\log |s| =[\widetilde{\Delta}_{l}] -dd_{\widetilde{J}}^{c}\log |s|+ d(d^{c}\log |s|)\circ (\widetilde{J}-i)
\end{equation}
est dans la classe de cohomologie duale à la classe d'homologie de $\widetilde{\Delta}_{l}$. Comme $\widetilde{J}=i$ le long de $\widetilde{\Delta}_{l}$, la $1$-forme $\beta = d^{c}(\log |s|)\circ (\widetilde{J}-i)$ est bornée.
En appliquant $T_{\widetilde{f}_{n}}$ à l'identité (\ref{Lelong-Poincaré}), nous obtenons, grâce aux formules de Jensen et de Stokes:
$$T_{\widetilde{f_{n}}}(\omega_{\widetilde{\Delta}_{l}})=I_{\widetilde{f_{n}}}(\widetilde{\Delta}_{l}) - \frac{1}{2\pi}\int_{0}^{2\pi}\log |s\circ f_{n}(e^{i\theta})|d\theta + \log |s\circ f_{n}(0)| + O(L_{\widetilde{f_{n}}}).$$
Comme $\widetilde{f}_{n}(0)$ reste loin de $\widetilde{\Delta}_{l}$, $\log |s\circ \widetilde{f_{n}}(0)|$ est minoré, ce qui conclut la démonstration. \qed
\end{demo}
\hspace{1cm}

Remarquons enfin, en vue du paragraphe 4, que pour $0\leq k\leq 6$, $[\widetilde{T}].[E_{S_{k}}]\geq 0$, puisque $\widetilde{J}$ coïncide avec $i$ le long des diviseurs exceptionnels.\\

\section{Le courant $T$ est diagonal}

Le courant $T$ admet la décomposition suivante: $T=T_{diag} + T'$, où $T_{diag}=\alpha_{1} [\Delta_{1}] + \alpha_{2} [\Delta_{2}] + \alpha_{3} [\Delta_{3}]$ est la composante diagonale singulière de $T$. L'objet de ce paragraphe est de montrer que $T'$ est nul. Raisonnons par l'absurde et supposons $T' \neq 0$. L'ingrédient principal de la démonstration est le lemme suivant, analogue pour le support de $T'$ du "lemme géométrique" de Julien Duval relatif à l'adhérence d'une courbe entière évitant cinq droites (c.f. \cite{Green presque complexe}):

\begin{Lemgeom}\label{geom}
Par tout point du support de $T'$ passe un $J$-disque non trivial contenu dans le support de $T'$ et obtenu comme limite d'une suite de $J$-disques $d_{n}\subset f_{n}(\mathbb{D})$. 
\end{Lemgeom}

Nous démontrons ce lemme en nous ramenant à la dimension $1$ à l'aide des projections centrales $\pi_{S_{k}}$.

\subsection{Notations}

Soit $P\in \mathbb{P}^{2}(\Cplx )$ un point appartenant au support de $T'$. Les droites $L_{k}$ étant en position générale, il existe au moins deux droites distinctes de la configuration $C$, disons $L_{1}$ et $L_{2}$, telles que $P\not \in L_{1}\cup L_{2}$. Notons $S_{1}$ leur point d'intersection, et $L_{3}$, $L_{4}$ les deux autres droites de $C$. Soient alors $\delta$, $\delta_{1}$ et $\delta_{2}$ des disques ouverts de $\mathbb{P}^{1}(\Cplx )$ d'adhérences deux à deux disjointes, contenant respectivement les points $\pi_{S_{1}}(P)$, $\pi_{S_{1}}(L_{1}\setminus \{S_{1}\})$ (que nous noterons $0$) et $\pi_{S_{1}}(L_{2}\setminus \{S_{1}\})$ (que nous noterons $\infty$).\\

Comme les $J$-disques $f_{n}(\mathbb{D})$ évitent les droites $L_{1}$ et $L_{2}$, pour tout $n$, $\varphi_{1,n}(\mathbb{D} )$ évite les points $0$ et $\infty$. Ainsi, d'après l'inégalité (\ref{eq2}), le nombre moyen de feuillets du disque $\varphi_{1,n}(\mathbb{D})$  au-dessus de $\mathbb{P}^{1}(\Cplx )$ est équivalent au nombre d'îles de degré $1$ au-dessus de $\delta$ quand $n$ tend vers l'infini. Plus précisémment, si l'on note $I_{n,t}$ l'ensemble des îles de degré $1$ au-dessus de $\delta$ pour la fonction $\varphi_{1,n}$ contenues dans le disque fermé $\overline{D(0,t)}$, on a alors:
\begin{equation}\label{feuillets-îles}
\lim_{n\rightarrow +\infty} \frac{T_{\varphi_{1,n}}(\omega ')}{\int_{0}^{1}\#I_{n,t}\frac{dt}{t}}=1.
\end{equation}
Si $d$ est une île de degré $1$ au-dessus de $\delta$, la restriction de la projection $\pi_{S_{1}}$ à $f_{n}(d)$ est un homéomorphisme sur $\delta$ (voir Figure 1).\\

Si $T_{\varphi_{1,n}}(\omega ')$ n'est pas négligeable devant $T_{f_{n}}(\omega )$, alors de telles îles sont nombreuses relativement à $T_{f_{n}}(\omega )$, et nous pourrons définir le disque recherché comme une limite de sous-disques de la famille  $(f_{n}(\mathbb{D}))$ de degré $1$ au-dessus de $\delta$. Sinon, nous verrons que le courant $T'$ est localement vertical relativement à $\pi_{S_{1}}$, et nous conclurons à l'aide d'une projection depuis un autre point double. \\

Par la suite nous noterons $T_{1}$ le courant suivant:
$$T_{1} = \lim_{n\rightarrow +\infty} \frac{1}{T_{f_{n}}(\omega )}\int_{0}^{1}\left ( \sum_{d\in I_{n,t}}[f_{n}(d)]  \right )\frac{dt}{t}.$$
Distinguons deux cas, selon que $P$ appartienne ou non au support de $T_{1}$.

\subsection{Cas où $P\in \operatorname{Supp}T_{1}$}

Nous allons construire le disque recherché comme limite de $J$-disques $(f_{n}(d_{n}))$, où les $d_{n}$ sont des îles de degré $1$ au-dessus de $\delta$. Notons:
$$\mathcal{G} = \{ J\text{ -disques } \gamma \in \mathbb{P}^{2}(\Cplx ) \text{ | } (\pi_{S_{1}})_{|\gamma}:\gamma \rightarrow \delta \text{ est de degré 1 et } \gamma \cap C = \emptyset \}.$$ 

\begin{Fait}
L'ensemble $\mathcal{G}$ est relativement compact pour la distance de Haussdorff, et son adhérence $\overline{\mathcal{G}}$ vaut:
$$\overline{\mathcal{G}}=\mathcal{G} \cup (\pi_{S_{1}}^{-1}(\delta ) \cap L_{3})\cup (\pi_{1}^{-1}(\delta ) \cap L_{4}) \cup \{ S_{1} \}.$$
\end{Fait}
En effet, par projection centrale en $S_{1}$ des $J$-disques de $\mathcal{G}$, nous obtenons une famille d'homéomorphismes $\alpha$-quasiconformes $(\varphi_{\gamma})_{\gamma \in \mathcal{G}}$ à valeurs dans $\delta$, donc normale (voir le paragraphé 2.2.1). Le paragraphe 2.2.2 montre alors que $\mathcal{G}$ est une famille normale. Par positivité d'intersecion, son adhérence s'obtient en ajoutant simplement les disques $\pi_{S_{1}}^{-1}(\delta ) \cap L_{3}$ et $\pi_{S_{1}}^{-1}(\delta ) \cap L_{4}$ ainsi que le point $ S_{1}$.\\

Montrons qu'il existe un disque holomorphe $\gamma \in \overline{\mathcal{G}} \setminus \{ S_{1} \}$ passant par $P$ et contenu dans le support de $T_{1}$. Pour cela, considérons l'ensemble $\mathcal{M}(\overline{\mathcal{G}}\setminus \{ S_{1} \} )$ des mesures localement finies sur $\overline{\mathcal{G}}\setminus \{ S_{1} \}$, qui est un espace compact pour la topologie de la convergence faible sur tout compact de $\overline{\mathcal{G}}\setminus \{ S_{1} \}$. Pour un disque $\gamma \in \overline{\mathcal{G}}\setminus \{S_{1} \}$, $\delta_{\gamma}$ désignera la masse de Dirac au point $\gamma$. Pour tout $n\geq 0$, la mesure $\mu_{n}$ définie par:
$$\mu_{n} = \frac{1}{T_{f_{n}}(\omega )}\int_{0}^{1}\left ( \sum_{d\in I_{n,t}} \delta_{f_{n}(d)} \right )\frac{dt}{t}$$
est bien localement finie, car pour tout compact $K$ de $\overline{\mathcal{G}}\setminus \{ S_{1} \}$, les $J$-disques $\gamma \in K$ sont d'aire minorée en vertu du théorème de Lelong pour les $J$-disques (voir l'article de J.C. Sikorav dans \cite{Audin}). La suite $(\mu_{n})$ converge donc, quitte à en extraire une sous-suite, vers une mesure $\mu$  de $\mathcal{M}(\overline{\mathcal{G}}\setminus \{ S_{1} \} )$. On vérifie alors aisément le fait suivant:

\begin{Faitt}
$T_{1}=\displaystyle{\int_{\mathcal{G} \setminus \{ S_{1} \}}}[\gamma ] d\mu (\gamma )$. En particulier, $\operatorname{Supp}T_{1}=\displaystyle{\bigcup_{\gamma \in \operatorname{Supp}\mu}}\gamma$.
\end{Faitt}
Ainsi, comme $P$ est dans le support de $T_{1}$, il existe un disque $\gamma$ appartenant au support de $\mu$ passant par $P$. Comme $\operatorname{Supp} \mu \subset \overline{\bigcup_{n} \operatorname{Supp} \mu_{n}}$, $\gamma$ est bien limite de $J$-disques de $(f_{n}(\mathbb{D}))$.

\subsection{Cas où $P\not \in \operatorname{Supp}T_{1}$}

Nous allons démontrer que le courant $T'$ est vertical relativement à la projection $\pi_{S_{1}}$. Commençons par vérifier que $T_{1}$ est localement la partie horizontale de $T$ relativement à $\pi_{S_{1}}$:
\begin{Faittt}\label{fait 3}
$(\pi_{S_{1}})_{*}(T_{|\pi_{S_{1}}^{-1}(\delta )})=(\pi_{S_{1}})_{*}T_{1}$.
\end{Faittt}

\begin{demo}
D'après (\ref{feuillets-îles}) et la proposition \ref{comparaison projection}, on a:
$$\frac{1}{T_{f_{n}}(\omega )}\left ( \int_{0}^{1} \int_{D(0,t) \setminus \bigcup_{d\in I_{n,t}}}(\varphi_{1,n})^{*}\omega ' \right ) \frac{dt}{t} 
 = \frac{o(T_{\varphi_{1,n}}(\omega '))}{T_{f_{n}}(\omega )} = o(1).$$
On en déduit:
\begin{align}
(\pi_{S_{1}})_{*}(T_{|\pi_{S_{1}}^{-1}(\delta )}) & = \lim_{n\rightarrow +\infty}\frac{1}{T_{f_{n}}(\omega )} \int_{0}^{1} [\varphi_{1,n}(D(0,t))] \frac{dt}{t} \notag \\
&\ = \lim_{n\rightarrow +\infty}\frac{1}{T_{f_{n}}(\omega )} \int_{0}^{1} \left (\sum_{d\in I_{n,t}} \int_{d} [\varphi_{1,n}(d)] \right ) \frac{dt}{t} \notag \\
&\ = (\pi_{S_{1}})_{*}T_{1}. \notag
\end{align}
\qed
\end{demo}

Ainsi, comme $P$ n'est pas dans le support de $T_{1}$, $\mathds{1}_{B(P,\varepsilon )}T_{1}=0$ pour un certain $\varepsilon >0$, et donc: $(\pi_{S_{1}})_{*}(\mathds{1}_{B(P,\varepsilon )}T')=(\pi_{S_{1}})_{*}(\mathds{1}_{B(P,\varepsilon )}T)=0$. Quitte à réduire $\varepsilon$, on peut trouver dans $B(P,\varepsilon )$ un repère local $(\zeta_{1}, \zeta_{2})$ du fibré $T_{J}^{1,0}$ tel que $\zeta_{1}$ soit tangent aux $J$-droites du pinceau en $S_{1}$. Dans le repère dual, $\mathds{1}_{B(P,\varepsilon )}T'$ s'écrit: 
$$\mathds{1}_{B(P,\varepsilon )}T'=iT_{1,1}\zeta_{1}^{*}\wedge \overline{\zeta}_{1}^{*} + iT_{1,2}\zeta_{1}^{*}\wedge \overline{\zeta}_{2}^{*} + iT_{2,1}\zeta_{2}^{*}\wedge \overline{\zeta}_{1}^{*} + iT_{2,2}\zeta_{2}^{*}\wedge \overline{\zeta}_{2}^{*}.$$

Comme $(\pi_{S_{1}})_{*}(\mathds{1}_{B(P,\varepsilon )}T')=0$, alors $T_{2,2}=0$, puis $T_{1,2} = T_{2,1}=0$ (voir paragraphe 2.3.2). Comme de plus $T'$ est fermé, la fonction $T_{1,1}$ est constante le long des $J$-droites du pinceau en $S_{1}$, donnant une mesure positive $dT_{1,1}$ sur $\delta$. Le courant $\mathds{1}_{B(P,\varepsilon )}T'$ est donc de la forme:
\begin{equation}\label{vertical}
\mathds{1}_{B(P,\varepsilon )}T' = \int_{\pi_{S_{k}}(B(P,\varepsilon ))} [\pi_{S_{k}}^{-1}(z)\cap B(P,\varepsilon ) ]dT_{1,1}(z).
\end{equation}

Achevons la démonstration du lemme géométrique. Supposons pour commencer que $P$ n'appartienne pas à $L_{3}$. Notons $S_{3}$ le point double $L_{1}\cap L_{3}$. Comme $(\pi_{S_{3}})_{*}\left ([\pi_{S_{1}}^{-1}(z) ] \right )$ est non nul pour tout $z \in \delta$, l'égalité (\ref{vertical}) implique que le point $P$ est dans le support du courant $(\pi_{S_{3}})_{*}T'$. Comme dans le paragraphe 4.1, nous pouvons définir, à partir des îles de degré $1$ au-dessus d'un petit disque contenant $\pi_{S_{3}}(P)$ pour l'application quasiconforme $\varphi_{3,n}$, un courant $T_{3}$ de support inclus dans celui de $T'$. Le fait 3 s'applique, et $P$ est dans le support de $T_{3}$: le $J$-disque $\gamma$ recherché est alors construit comme dans le paragraphe 4.2. La situation est identique si $P$ n'appartient pas à $L_{4}$.\\

Il reste à traiter le cas où $P$ est le point double $S_{2}=L_{3}\cap L_{4}$ situé en vis-à-vis de $S_{1}$. Nous allons voir que ce point appartient forcément au support de $T_{1}$, et le paragraphe 4.2 permettra alors de conclure la démonstration du lemme géométrique.\\
Supposons que ce ne soit pas le cas; l'équation (\ref{vertical}) est alors vérifiée. La diagonale $\Delta_{1}$ reliant les points $S_{1}$ et $S_{2}$ n'est pas une composante singulière de $T'$, donc le point $\pi_{S_{1}}(S_{2})$ ne peut être un atome pour la mesure $dT_{1,1}$. Soit $z\in \operatorname{Supp}(dT_{1,1})\setminus \{\pi_{S_{1}}(S_{2}) \}$ un point proche de $\pi_{S_{1}}(S_{2})$.\\ 
Comme nous l'avons vu plus haut, $(\pi_{S_{3}})_{*}[\overline{\pi_{S_{1}}^{-1}(z)}]\neq 0$, donc le point $P_{z}=\pi_{S_{1}}^{-1}(z)\cap L_{4}$ est dans le support de $T_{3}$: il existe un $J$-disque $\gamma_{z}$ contenu dans le support de $T_{3}$, issu de la suite $(f_{n}(\mathbb{D}))$, passant par $P_{z}$. Par positivité d'intersection, $\gamma_{z}\subset L_{4}$, donc comme $(\pi_{S_{1}})_{*}[L_{4}] \neq 0$, le courant $T_{1}$ a de la masse près de $P_{z}$: le point $z$ pouvant être pris arbitrairement proche de $\pi_{S_{1}}(S_{2})$, le point $S_{2}$ est dans le support de $T_{1}$, contrairement à ce qui a été supposé. \\

\subsection{Conclusion à l'aide du lemme géométrique}

Rappelons que nous avons supposé $T'$ non nul. Nous allons obtenir une contradiction à l'aide du lemme géométrique et du fait suivant, analogue du théorème de  Liouville pour le courant $T'$:

\begin{Faitttt}
Toute $J$-droite de $\mathbb{P}^{2}(\Cplx )$ rencontre le support de $T'$.
\end{Faitttt}

\begin{demo}
Soit $L_{0}$ une $J$-droite de $\mathbb{P}^{2}(\Cplx )$ non entièrement contenue dans le support de $T'$. Notons $(L_{z})_{z\in \mathbb{P}^{1}(\Cplx )}$ le pinceau de $J$-droites centré en un point de $L_{0}\setminus \operatorname{Supp}(T')$. Le support de $T'$ rencontre forcément une $J$-droite $L_{z}$ du pinceau, et par tout point de $L_{z}\cap \operatorname{Supp}(T')$ passe un $J$-disque contenu dans $\operatorname{Supp}(T')$. Au moins un  $J$-disque $\gamma \subset \operatorname{Supp}(T')$ rencontre $L_{z}$ sans y être contenu (sinon, le fait pour un point de $L_{z}$ d'être contenu dans un $J$-disque inclus dans $\operatorname{Supp}(T')\cap L_{z}$ serait une propriété ouverte, et $L_{z}$ serait contenue dans le support de $T'$, ainsi que le centre du pinceau). Par positivité d'intersection, les droites proches de $L_{z}$ coupent aussi $\gamma$, et finalement toutes les droites du pinceau, et en particulier $L_{0}$, rencontrent le support de $T'$.
\qed 
\end{demo}

L'ensemble $L_{1} \cap \operatorname{Supp}T'$ est donc un ensemble fermé non vide. C'est aussi un ouvert de $L_{1}$, car par tout point $P$ de $L_{1} \cap \operatorname{Supp}T'$ passe un disque holomorphe $\gamma$ non trivial contenu dans $\operatorname{Supp}T'$ et issu de $(f_{n}(\mathbb{D}))$, donc contenu dans $L_{1}$ par positivité d'intersection. Ainsi, $L_{1} \subset \operatorname{Supp}T'$, et de même $L_{2} \subset \operatorname{Supp}T'$. Par le point $L_{1}\cap L_{2}$ passe donc un $J$-disque non trivial issu de $(f_{n}(\mathbb{D})$. Toujours par positivité d'intersection, ce $J$-disque doit être contenu à la fois dans $L_{1}$ et dans $L_{2}$, ce qui est impossible.

\section{Obtention d'une contradiction aux points doubles}

Le paragraphe précédent montre que le courant $T$ est porté par $\Delta$. Suivant une idée de M.L. McQuillan (c.f. \cite{Bloch hyperbolicity}), nous allons examiner le comportement des $J$-disques $f_{n}$ au voisinage des points doubles de la configuration $C$ à l'aide d'éclatements. Rappelons que $(\widetilde{X}, \widetilde{J})$ désigne l'éclaté du plan projectif presque complexe $(\mathbb{P}^{2}(\Cplx ),J)$ en les $S_{k}$, que les $E_{S_{k}}$ ($1\leq k \leq 6$) désignent les diviseurs exceptionnels au-dessus des points $S_{k}$, et que $\Pi:\widetilde{X}\rightarrow \mathbb{P}^{2}(\Cplx )$ est la projection canonique. Les longueurs et les aires seront calculées pour la métrique hermitienne associée à la $(1,1)$-forme $\widetilde{\omega}$ définie au paragraphe 2.1.3.\\

Le courant $\widetilde{T}$ est porté par la réunion de $\widetilde{J}$-courbes $D=(\bigcup \widetilde{\Delta_{l}})\cup (\bigcup E_{S_{k}})$. Comme il est fermé et d'ordre zéro, il est de la forme: 
\begin{equation}\label{forme}
\widetilde{T}= \widetilde{\alpha_{1}}[\widetilde{\Delta_{1}}] + \widetilde{\alpha_{2}}[\widetilde{\Delta_{2}}] + \widetilde{\alpha_{3}}[\widetilde{\Delta_{3}}] + \sum_{1\leq k \leq 6}\lambda_{k}[E_{S_{k}}].
\end{equation}
Pour tout $1\leq k\leq 6$, $[\widetilde{T}].[E_{S_{k}}]\geq 0$ (voir paragraphe 3.2), donc $\widetilde{T}$ ne peut être porté par les diviseurs exceptionnels, qui sont d'autointersection $-1$. \\
Rappelons que $S_{1}$ est le point d'intersection de $L_{1}$ et $L_{2}$. La forme particulière de $\widetilde{T}$ (voir \ref{forme}) permet de calculer $[\widetilde{T}].[\widetilde{\Delta}_{1}]$:
$$[\widetilde{T}].[\widetilde{\Delta}_{1}]=-\widetilde{\alpha}_{1} + \widetilde{\alpha}_{2} + \widetilde{\alpha}_{3} + \lambda_{1} + \lambda_{2} \geq I(\widetilde{T}, \widetilde{\Delta}_{1}).$$
La contradiction recherchée découle du lemme suivant:

\begin{Lem}
$\widetilde{\alpha}_{2} + \widetilde{\alpha}_{3}+ \lambda_{1} + \lambda_{2} \leq I(\widetilde{T}, \widetilde{\Delta}_{1}) $.
\end{Lem}

\begin{demo}
Nous allons établir que $\lambda_{1}$ est majoré par l'intersection géométrique asymptotique près de $E_{S_{1}}$ du courant $\widetilde{T}$ et de $\widetilde{\Delta}_{1}$ (le raisonnement se transposera aisément aux cas des diviseurs $E_{S_{2}}$, $\widetilde{\Delta}_{2}$ et $\widetilde{\Delta}_{3}$). Nous allons pour cela projeter sur $E_{S_{1}}$ les morceaux de $\widetilde{f}_{n}(\mathbb{D})$ situés près du diviseur exceptionnel, de façon à pouvoir appliquer la théorie d'Ahlfors. L'inégalité recherchée découlera alors de l'inégalité \ref{eq1} et du fait que les $\widetilde{J}$-disques $\widetilde{f}_{n}$ évitent $\widetilde{L}_{1}$ et $\widetilde{L}_{2}$: le nombre d'îles au-dessus du point $\widetilde{\Delta}_{1}\cap E_{s_{1}}$ majore asymptotiquement la masse près de $E_{S_{1}}$, soit environ, après intégration à la Nevanlinna, $\lambda_{1}T_{\widetilde{f}_{n}}(\widetilde{\omega})$.\\

Commençons par localiser l'étude au voisinage de $E_{S_{1}}$. Il existe une famille $(U_{t})_{0<t<\eta}$ de voisinages pseudoconvexes de $E_{S_{1}}$, le voisinage $U_{\eta}$ étant pris suffisamment petit. En identifiant $E_{S_{1}}$ avec $\mathbb{P}^{1}(\Cplx )$, la projection qui a tout point $\widetilde{P}$ de $U_{\eta}$ associe le point d'intersection de la transformée stricte de la $J$-droite passant par $\Pi (\widetilde{P})$ et $S_{1}$ avec $E_{S_{1}}$ s'indentifie à la projection $\widetilde{\pi}_{S_{1}}$. \\
Les longueurs et les aires sur le diviseur exceptionnel $E_{S_{1}}$ seront calculées pour la métrique asociée à la forme de Fubini-Study $\omega '$.\\ 

Fixons un $0<\varepsilon <\eta$. Posons $\Sigma_{n} = \widetilde{f_{n}}^{-1}(U_{\varepsilon})$. Les applications $\widetilde{\pi}_{E_{S_{1}}} \circ \widetilde{f}_{n} : \Sigma_{n} \rightarrow E_{S_{1}}$ s'identifient aux applications $\varphi_{1,n}$ et sont $\alpha$-quasiconformes. Posons, pour tout $0<t<1$, $\Sigma_{n,t} = \Sigma_{n} \cap D(0,t)$. Nous allons, grâce à la théorie d'Ahlfors, établir l'inégalité:
\begin{equation}\label{fin}
\int_{0}^{1}\operatorname{Aire}(\varphi_{1,n}(\Sigma_{n,t}))\frac{dt}{t} + o(T_{\widetilde{f}_{n}}(\widetilde{\omega})) \leq \int_{0}^{1}\#[\widetilde{f}_{n}(D(0,t))].[\widetilde{\Delta}_{1} \cap U_{\varepsilon}]\frac{dt}{t}. 
\end{equation}
Comme $(\widetilde{\pi}_{S_{1}})_{*}(\mathds{1}_{U_{\varepsilon}}\widetilde{T})=\lambda_{1}[E_{1}]$, nous obtiendrons, après division par $T_{\widetilde{f}_{n}}(\widetilde{\omega})$:
$$\lambda_{1} \leq \lim \sup_{n\rightarrow +\infty} \frac{1}{T_{\widetilde{f}_{n}}(\widetilde{\omega})} \int_{0}^{1}\#[\widetilde{f}_{n}(D(0,t))].[\widetilde{\Delta}_{1} \cap U_{\varepsilon}]\frac{dt}{t}$$
ce qui est l'inégalité recherchée.\\

Supposons dans un permier temps que $\Sigma_{n,t}$ soit connexe. Commençons par montrer que \\
$\int_{0}^{1}-\min (0, \chi (\Sigma_{n,t}))\frac{dt}{t}$ est faible devant $T_{\widetilde{f}_{n}}(\widetilde{\omega})$. Le domaine $\Sigma_{n,t}$ est de la forme $\Omega_{n} \setminus \bigcup_{j \leq j_{n,t}}\overline{\Omega_{n,j}}$, où $j_{n,t}$ est un entier positif, $\Omega_{n}$ ainsi que les $\Omega_{n,j}$ sont des disques topologiques ouverts, et où pour tout $j$: $\widetilde{f}_{n}(\overline{\Omega_{n,j}})\cap \overline{U_{\varepsilon}}=\widetilde{f}_{n}(\partial \Omega_{n,j})$ (voir Figure 2).\\

Comme nous l'avons vu au paragraphe 2.3.2, pour tout $1\leq j \leq j_{n,t}$, le disque $\Omega_{n,j}$ ne peut avoir son image $\widetilde{f}_{n}(\Omega_{n,j})$ entièrement contenue dans $U_{\eta}$, puisque les ouverts $U_{t }$ sont tous pseudo-convexes. Le $\widetilde{J}$-disque $\widetilde{f}_{n}(\Omega_{n,j})$ ne peut non plus rester près du diviseur $D$. En effet, si $\widetilde{f}_{n}(\Omega_{n,j})$ était contenu dans le $2\varepsilon$-voisinage $(D)_{2\varepsilon}$ de $D$, l'image du $J$-disque $f_{n}(\Omega_{n,j})$ par une projection centrale sur $\Delta_{1}$ depuis un point $P$ situé hors de $\Delta_{1}$ contiendrait forcément le point $S_{2}$ , en vertu du principe du maximum pour les applications quasiconformes. Le $J$-disque $f_{n}(\Omega_{n,j})$ passerait donc près du diviseur exceptionnel $E_{S_{2}}$, coupant ainsi les droites $L_{3}$ et $L_{4}$, ce qui est impossible.\\

Ainsi, pour tout $j$, le $\widetilde{J}$-disque $\widetilde{f}_{n}(\Omega_{n,j})$ sort de $(D)_{2\varepsilon}$, et d'après le théorème de Lelong, il existe une constante $c_{\varepsilon}>0$ telle que l'aire de la portion de $\widetilde{f}_{n}(\Omega_{n,j})$ contenue dans $(D)_{2\varepsilon} \setminus (D)_{\varepsilon}$ soit minorée par $c_{\varepsilon}$. Comme $\widetilde{T}$ ne charge pas $(D)_{2\varepsilon} \setminus (D)_{\varepsilon}$, de tels disques sont peu nombreux et $\int_{0}^{1}-\min (0, \chi (\Sigma_{n,t}))\frac{dt}{t}$ est bien négligeable devant $T_{f_{n}}(\widetilde{\omega})$.\\

Vérifions également que la longueur de l'image du bord de $\Sigma_{n,t}$ est négligeable devant le terme d'aire. Plus précisément, nous n'aurons besoin de nous intéresser qu'à la partie de $\partial \Sigma_{n,t}$ dont l'image par $\widetilde{f}_{n}$ est loin de $\widetilde{\Delta_{1}}$: posons $C_{n,t} = \partial \Sigma_{n,t} \setminus \widetilde{f}_{n}^{-1}\left ((\widetilde{\Delta}_{1})_{\varepsilon}\right )$, $(\widetilde{\Delta}_{1})_{\varepsilon}$ désignant le $\varepsilon$-voisinage de $\widetilde{\Delta}_{1}$. D'une part, rappelons que $L_{\widetilde{f}_{n}}$ est négligeable devant $T_{\widetilde{f}_{n}}(\widetilde{\omega})$. D'autre part, comme $\widetilde{T}$ ne charge pas $U_{\eta} \setminus U_{\varepsilon}\cup (\widetilde{\Delta}_{1})_{\varepsilon}$, la formule de coaire permet d'affirmer, quitte à modifier légèrement $\varepsilon$, que $\int_{0}^{1}\operatorname{Longueur}(\widetilde{f}_{n}(C_{n,t}))\frac{dt}{t}$ est négligeable devant $T_{\widetilde{f}_{n}}(\widetilde{\omega})$.\\

En vue d'appliquer la théorie d'Ahlfors, considérons le disque ouvert $\delta = \widetilde{\pi}_{S_{1}}((\widetilde{\Delta}_{1})_{\varepsilon}\cap U_{\varepsilon})$ de $E_{S_{1}}$, contenant le point $\widetilde{\pi}_{S_{1}}(\widetilde{\Delta}_{1}\cap U_{\varepsilon})$, et posons $\Sigma_{n,t}' = \Sigma_{n,t} \setminus \varphi_{1,n}^{-1}(\overline{\delta})$.  Comme les disques $\widetilde{f}_{n}(\mathbb{D} )$ évitent les droites $\widetilde{L}_{1}$ et $\widetilde{L}_{2}$, les fonctions $\varphi_{1,n}$ évitent deux points distincts, que nous noterons $0$ et $\infty$. Nous disposons d'un recouvrement quasiconforme $\varphi_{1,n}:\Sigma_{n,t} \rightarrow E_{S_{1}}\setminus (\overline{\delta} \cup \{ 0, \infty \})$. L'inégalité (\ref{eq1}) est vérifiée: il existe une constante $h>0$ telle que si $I(n,t)$ désigne le nombre d'îles au-dessus de $\delta$ dans $\Sigma_{n,t}$:
\begin{equation}\label{der}
\operatorname{Aire}(\varphi_{1,n}(\Sigma_{n,t})) \leq I(n,t) -\min (0, \chi (\Sigma_{n,t})) + h\operatorname{Longueur}(\varphi_{1,n}(\partial \Sigma_{n,t}' )\setminus \partial \delta ).
\end{equation}

Comme d'une part, l'image par $\widetilde{f}_{n}$ de toute île au-dessus de $\delta$ coupe $\widetilde{\Delta}_{1}$, et d'autre part, $\operatorname{Longueur}(\varphi_{1,n}(\partial \Sigma_{n,t}' )\setminus \partial \delta )$ est majorée (à une contante multiplicative près dépendant de $\widetilde{\omega}$) par $\operatorname{Longueur}(\widetilde{f}_{n}(C_{n,t}))$, on obtient l'inégalité:

$$\operatorname{Aire}(\varphi_{1,n}(\Sigma_{n,t})) \leq \#[\widetilde{f}_{n}(D(0,t))].[\widetilde{\Delta}_{1} \cap U_{\varepsilon}]  - \min (0,\chi (\Sigma_{n,t})) + \operatorname{Longueur}(\widetilde{f}_{n}(C_{n,t}))$$
qui, après intégration par rapport à $\frac{dt}{t}$, donne l'inégalité (\ref{fin}) recherchée.\\

Si l'on ne suppose plus $\Sigma_{n,t}$ connexe, il existe une inégalité du type (\ref{der}) pour chacune de ses composantes connexes; il suffit alors de sommer ces inégalités.\\

Nous pouvons reproduire ce raisonnement à l'identique au voisinage du diviseur $E_{2}$, puisque les disques $\widetilde{f}_{n}$ évitent $\widetilde{L}_{3}$ et $\widetilde{L}_{4}$, ainsi qu'au voisinage des diviseurs  $\widetilde{\Delta}_{2}$ et $\widetilde{\Delta}_{3}$, puisque les diviseurs exceptionnels $E_{3}$, $E_{4}$, $E_{5}$ et $E_{6}$ sont aussi évités. Ceci achève la démonstration du lemme, puis du théorème de Bloch. \qed
\end{demo}
\vspace{1cm}
\begin{center}
  \large\bf Figures
\end{center}

\begin{picture}(200,85)(-45,-5)
\linethickness{0.05mm}
\put(-5,15){\line(1,0){70}}
\put(30,35){\circle*{0.8}}
\multiput(30,65)(-1,-3){20}{\line(-1,-3){0.5}}
\multiput(30,65)(1,-3){20}{\line(1,-3){0.5}}
\put(26,5){$(S_{1}P)$}
\put(26,-5){\bf Figure 1}
\put(31,34){$P$}
\put(31.5,63.5){$S_{1}$}
\put(35,75){$L_{1}$}
\put(22,75){$L_{2}$}
\put(-5,38){$L_{3}$}
\put(71,40){$L_{4}$}
\put(0,12){$0$}
\put(56.5,13){$\infty$}
\put(66,13){$\mathbb{P}^{1}(\Cplx )$}
\put(48,25){$\pi_{S_{1}}$}
\put(31,12){$\delta$}

\linethickness{0.30mm}
\put(30,10){\line(0,1){63.5}}
\put(13.3,15){\line(1,0){33.4}}
\put(13.3,14){\line(0,1){2}}
\put(46.7,14){\line(0,1){2}}
\put(30,65){\line(-3,-5){32}}
\put(30,65){\line(3,5){5}}
\put(30,65){\line(-3,5){5}}
\put(0,38){\line(7,-1){70}}
\put(70,40){\line(-5,-1){70}}
\put(30,65){\line(3,-5){32}}
\cbezier(15,18.38)(25,16)(37,22)(44.46,20)
\cbezier(17,24.36)(25,22)(35,27)(43.46,26)
\cbezier(25,48.28)(27,45)(32,52)(34.43,50)
\cbezier(23,42.3)(29,40)(31,48)(36.10,45)
\cbezier(22,39.31)(28,36)(32,44)(37.11,42)
\put(45,30){\vector(3,-7){6}}

\end{picture}

\vspace{2cm}

\begin{picture}(200,60)(-10,-5)
\linethickness{0.5mm}
\put(132,24){$E_{S_{1}}$}
\put(51,25){$U_{\varepsilon}$}
\put(97,-2.5){$(\widetilde{\Delta}_{1})_{\varepsilon}$}
\put(55,-3){\bf Figure 2}
\put(99,46){$\widetilde{\Delta}_{1}$}
\put(73,27){$\widetilde{\pi}_{S_{1}}$}
\put(-3,34){$\Omega_{n,1}$}
\put(14,24){$\Omega_{n,2}$}
\put(7,43){$D(0,t)$}
\put(74,49){$\widetilde{f}_{n}(\Omega_{n,1})$}
\put(110,39){$\widetilde{f}_{n}(\Omega_{n,2})$}
\put(37,26.5){$\widetilde{f}_{n}$}
\linethickness{0.05mm}
\put(100,5){\line(0,1){40}}
\put(60,25){\line(1,0){70}}
\multiput(60,30)(1,0){70}{\line(1,0){0.5}}
\multiput(60,20)(1,0){70}{\line(1,0){0.5}}
\multiput(95,5)(0,1){40}{\line(0,1){0.5}}
\multiput(105,5)(0,1){40}{\line(0,1){0.5}}
\put(80,29){\vector(0,-1){4}}
\linethickness{0.30mm}
\cbezier(59,30)(57.5,29)(59.80,26)(56.5,25)
\cbezier(59,20)(57.5,21)(59.80,24)(56.5,25)
\cbezier(95,4)(96,2.5)(99,4.80)(100,1.5)
\cbezier(105,4)(104,2.5)(101,4.80)(100,1.5)
\cbezier(72,30)(77,70)(87,27)(95,27)
\cbezier(95,27)(100,26)(110,50)(120,28)
\cbezier(72,30)(71,25)(67,26)(60,28)
\cbezier(120,28)(121,25.5)(122,25)(130,26)
\put(10,30){\circle{35}}
\put(0,35){\circle{10}}
\put(17,25){\circle{12}}
\put(30,25){\vector(1,0){15}}
\end{picture}


\begin{thebibliography}{2}
   \bibitem{Audin} \textsc{M. Audin et J. Lafontaine ed.}, \emph{Holomorphic curves in symplectic geometry}, Progress in Math., vol 117, Birkhäuser, Basel, 1994.
   \bibitem{Brody} \textsc{R. Brody}, \emph{Compact manifolds and hyperbolicity}, Trans. Amer. Math. Soc. 235 (1978), 213-219.
   \bibitem{Brunella} \textsc{M. Brunella}, \emph{Courbes entières et feuilletages holomorphes}, L'enseignement mathématique (1999), no. 45, 195-216.
   \bibitem{Green presque complexe} \textsc{J. Duval}, \emph{Un théorème de Green presque complexe}, Annales de l'institut Fourier 54 no 7 (2004), 2357-2367.
   \bibitem{Elkhadhra} \textsc{F. Elkhadhra} \emph{J-pluripolar subsets and currents on almost complex manifolds}, Math. Zeit. (2010) vol. 264 no. 2, 399-422.
   \bibitem{Green} \textsc{M. Green}, \emph{Some Picard theorems for holomorphic maps to algebraic varieties}, Amer. J. Math. 97 1975), 43-75.
   \bibitem{Gromov} \textsc{M. Gromov}, \emph{Pseudo holomorphic curves in symplectic manifolds}, Ivent.Math.83 (1985), 3007-347.
   \bibitem{Kruglikov} \textsc{B. Kruglikov et M. Overholt}, \emph{Pseudoholomorphic mappings and Kobayashi hyperbolicity}, Diff. Geom. Appl. 11 (1999), 265-277.
   \bibitem{Lehto} \textsc{O. Lehto and K.I. Virtanen}, \emph{Quasiconformal mappings in the plane}, Grund. der. math. Wiss, Vol. 126, Springer, Berlin, 1973.
   \bibitem{Bloch hyperbolicity} \textsc{M. McQuillan}, \emph{Bloch hyperbolicity}, preprint IHES (2001).
   \bibitem{Pali} \textsc{Nefton Pali}, \emph{Fonctions plurisousharmoniques et courants de type (1,1) sur les variétés presque complexes}, arXiv
   \bibitem{Sikorav} \textsc{J.-C. Sikorav}, \emph{Dual elliptic Planes}, preprint 2000 arXiv math.SG/008234.
\end{thebibliography}
\end{document}